\documentclass[12pt]{article}
\usepackage{amssymb,latexsym}
\begin{document}
\title{Bounding the torsion in CM elliptic curves}
\author{Dipendra Prasad and C S Yogananda}

\maketitle
\newtheorem{thm}{Theorem}[section]
\newtheorem{propose}{Proposition}[section]
\newtheorem{lemma}[propose]{Lemma}
\newenvironment{proof}{
                       \trivlist \item[\hskip \labelsep{\bf Proof}:]
                      }{
                        \hfill$\Box$\endtrivlist
                      }

\tableofcontents

{\sf

\begin{section}
{Introduction}

In \cite{as:cm} 
A. Silverberg, using the main theorem of complex multiplication of Shimura 
and Taniyama, has obtained a bound for the order of a 
point of finite order in a  CM abelian variety over a number field
 in terms of 
only the degree of the number field and the dimension of the abelian variety. 
As a corollary she obtains the following result for elliptic curves: {\it 
Let $E$ be an elliptic curve over a number field $K$ of degree $d$ with CM by 
an order $\mathfrak{O}$ in an imaginary quadratic field $k$. Suppose $P\in 
E(K)$ is a point of order $N$. Then $\varphi (N)\leq \delta \mu d$ where $\mu$
 is the number of roots 
of unity in $\mathfrak{O}$ and $\delta=1/2\mbox{ or }1$ depending on whether 
$k$ is contained in $K$ or not}.
Her results also imply a bound on the  full torsion subgroup 
of CM elliptic curves. 

The aim of this paper is to give an estimate on the order of the  torsion
subgroup of a CM  elliptic curve over a number field using only the result of
Deuring about supersingular primes, and elementary algebraic number theory.
To state our theorem, we need the following notation: if $M= l_1^{j_1}\cdots
l_r^{j_r}$ is the prime factorisation of $M $ let
$M^{\prime}=l_1^{(j_1+\delta_1)/2}\cdots l_r^{(j_r+\delta_r)/2}$ where
$\delta_i=0$ if $j_i$ is even and $\delta_i=1$ if $j_i$ is odd. 

 \begin{thm} Let $E$ be an elliptic curve over a number field $K$ of
degree $d$ with CM by an order in an imaginary quadratic field 
$k$.
 Then if $M$ is the order of the
 torsion subgroup of $E(K)$ we have:  \begin{enumerate} \item
$\varphi(M)\leq 2d$ if $K \cap k = {\bf Q} $;  \item $\varphi (M^{\prime}) 
\leq 2d$ if $k \subseteq K$ but $k\not = {\bf Q}(i),{\bf Q}(\omega ) $, 
($\omega$ being a third root of unity). 
\end{enumerate} \end{thm}

When $k = {\bf Q}(i),{\bf Q}(\omega )$ our method gives an extra factor of 
$2^{d(M)+1}$ where $d(M)$ is  the number of distinct prime divisors of $M$. 

\vspace{2mm}
\noindent{\bf Remark 1 :} A. Silverberg has remarked that one can easily 
give an estimate to the order of the full torsion subgroup of a CM elliptic
curve which is the same as our case $(i)$ from her theorem in [9]. 
Her estimate is better than ours in case $(ii)$. In case $(i)$, 
this follows
as the torsion subgroup on an elliptic curve is always of the form ${\bf  Z}/a
\times {\bf Z}/b$ where $a | b$. If the torsion subgroup of $E(K)$ 
contains ${\bf Z}/r \times {\bf Z}/r$ for $r> 2$,  it follows from
[10] that the field of complex multiplication (= $k$) must be contained 
in $K$. Therefore if $K \cap k = {\bf Q}$, then the torsion subgroup
is either ${\bf Z}/N$ or ${\bf Z}/N \times {\bf Z}/2$, where $N$ is the
maximum order of a torsion point, achieving the same bound that we obtain 
 for the order
of the torsion subgroup of $E(K)$ from the maximal order of a torsion element.

\vspace{2mm}
  
There is a large amount of literature on the torsion subgroup of elliptic 
curves over number fields. 
In \cite{lm:uni} Merel has shown that the order of the torsion subgroup of an 
elliptic curve over a number field $K$ can be bounded in terms of only the 
degree, $d$, of $K$ over {\bf Q}. The bound thus obtained (first by Merel and 
then improved by Oesterl\'e ) is exponential in $d$. The initial motivation
for this note was 
to investigate as to what could be the `right bound' by looking at the CM case
when we discovered that Silverberg has already done this.  
We also refer to
the paper of Olson [8] which deals with elliptic curves with complex
multiplication. 

\vspace{5mm}

\noindent {\bf Acknowledgement:} We are grateful  to Prof. J. Oesterl\'e whose
wonderful lectures on Merel's work and his own refinement at the Mehta 
Research Institute and at the ICTP, Trieste, stimulated our interest in 
working out the CM case. We are  grateful to Prof. J. Alperin for 
providing us with the proof of Lemma 2.3 and to Prof. A. Silverberg for 
remark 1. This work was done while the second
author was visiting Mehta Research Institute whose hospitality he
gratefully acknowledges.  \end{section}

\begin{section}
{Preliminary Lemmas}

\begin{lemma} Let $k$ be a quadratic extension of {\bf Q} and $K$ an
extension of {\bf Q} of degree $d$ with $K\cap k={\bf Q}$. Then the set of
primes $p$ in {\bf Q} which remain inert in $k$ and have the property that
there is at least one prime of degree 1 in $K$ above $p$ is of density
at least $1/(2d)$.  \end{lemma}

\begin{proof} Let $L$ be a Galois extension of {\bf Q} containing $K$ and
$k$, and let $G={\rm Gal}(L/{\bf Q})$. Further let $H_K$ and $H_k$ be the
subgroups of $G$ corresponding to the subfields $K$ and $k$, respectively,
of $L$. It is easy to see that the set of prime
 ideals $\mathfrak{p}$ in $L$, such that the prime ideal $\mathfrak{p}\cap
K$ is of degree 1 are precisely those for which the corresponding
Frobenius element $\sigma$ in $G$ belongs to $H_K$. The prime
$p=\mathfrak{p}\cap{\bf Q}$ is inert in $k$ if and only if $\sigma$ does
not belong to $H_k$. So, the primes $p$ in {\bf Q} as desired in the lemma
are precisely those for which there is a
prime $\mathfrak{p}$ in $L$ above $p$ for which the Frobenius 
element belongs to $(G\setminus
H_k)\cap H_K$. Since the cardinality of $(G\setminus H_k)\cap H_K$ is
$|H_K|/2$, it follows from the Cebotarev density theorem that the density
of $p$ in {\bf Q} as desired in the lemma is at least $1/(2d)$.  \end{proof}

\begin{lemma} Let $K$ be a number field of degree $d$ containing an
imaginary quadratic field $k$. Then the set of primes $p$ in {\bf Q} which
are inert in $k$ and have a prime of degree 2 in $K$ over $p$ is of
density at least $1/d$.  \end{lemma}

\begin{proof} Let $L$ be a Galois extension of {\bf Q} containing $K$ and
$G={\rm Gal}(L/{\bf Q})$. Further let $H_K$ and $H_k$ be the subgroups of
$G$ corresponding to the subfields $K$ and $k$, respectively, of $L$. The
set of primes $p$ as desired in the lemma are precisely those for which
the corresponding Frobenius substitution $\sigma$ does not belong to $H_k$
but whose square is in $H_K$.
Since $k$ is imaginary quadratic, the complex conjugation does not 
belong to $ H_k$. The following lemma combined with the Cebotarev density 
theorem completes the proof of our lemma. 
\end{proof}

\begin{lemma} Let $G$ be a finite group,  $N$ a subgroup of $G$
 of index 2 in $G$ and $H$ a subgroup of $N$. Suppose that there is an element 
of order 2, say $c$, in $G$ which is not in $N$. Then the set of elements 
in $G$ which do not belong to $N$ but whose square belongs to $H$ has 
cardinality at least that of $H$.
\end{lemma}

\begin{proof} (due to J. Alperin) We need to count elements $n\cdot c$ with
$n \in N$ whose square belongs to $H$. Clearly $HcH$ is a subset of $NcN = Nc$. We will prove that there are exactly $|H|$ elements in $HcH$ whose 
square belongs to $H$ which will prove our lemma. To prove this let $
A = H \cap cHc^{-1}$, and let $X \subset H$ be a set of left coset 
representatives of $A$ in $H$ so that every element of $H$ can be written 
uniquely in the form $x \cdot a$ with $x \in X$ and $a \in A$. From this it 
is easy to see that an element of $HcH$ can be uniquely written in the 
form $xch$ with $x \in X, h \in H$. Now $(xch)^2 = xchxch$ belongs to $H$ if and only if $chxc$ belongs to $H$ which happens if and only if $hx$ 
belongs to $cHc^{-1}$. Since both $x$ and $h$ belongs to $H$, we find that
$(xch)^2$ belongs to $H$ if and only if $hx$ belongs to $A=H \cap cHc^{-1}$.
For each $x$, this means that $h$ belongs to $Ax^{-1}$. So for each $x$,
there are $|A|$ many choices for $h$ such that $(xch)^2$ belongs to $H$.
Therefore the total number of elements in $HcH$ whose square belongs to $H$ is $$|A|\cdot \frac{|H|}{|A|} = |H|,$$ proving the lemma.     
\end{proof}

\end{section}

\begin{section}
{Proof of the main theorem}

 The proof of our main theorem will be a simple consequence of the lemmas 
in the previous section, Cebotarev density theorem, and  the
well known theorem about elliptic curves with complex multiplication 
that a prime $\mathfrak{p}$ in $K$ which is a prime of good
reduction for $E$ over $K$ is a prime of supersingular reduction if and
only if $\mathfrak{p}\cap {\bf Q}=p$ is inert or ramified in $k$ (see
\cite{sl:deu}).

\medskip

\noindent{\bf Case 1:} $K\cap k={\bf Q}$. 
We consider the set of rational primes $p$
coprime to $M$ which are inert in $k$ and have a split factor in $K$,
i.e., there is a prime of degree 1, say $\mathfrak{p}$ in $K$ which
divides $p$. Denote by ${\bf F}_{\mathfrak{p}}$
 the residue field associated to the prime ideal $\mathfrak{p}$ of $K$.
Then the torsion subgroup of $E$ over $K$ will inject into $E_{\mathfrak{p}}({\bf
F}_{\mathfrak{p}})$ which has cardinality $p+1$ since $\mathfrak{p}$ is
a prime of supersingular reduction of $E$ and ${\bf F}_{\mathfrak{ p}}$ is a finite field with $p$
elements; hence $M\vert (p+1)$. By Lemma 2.1 the density of such primes $p$
is at least $1/(2d)$ whereas by Cebotarev density theorem, the density of primes 
$p$ which are congruent to $-1$ modulo $M$ is $1/\varphi(M)$. Therefore 
we must have
$1/(2d)\leq 1/\varphi(M)$ and so $\varphi (M)\leq 2d$. 

\bigskip

\noindent {\bf Case 2:}  $k\subset K$.
We consider
primes $p$ in {\bf Q} such that $p$ is inert in $k$ and has a prime
factor $\mathfrak{p}$ of degree 2 in $K$ which is a prime of 
good reduction of $E$. The elliptic curve will have supersingular 
reduction at such primes. The density of such primes is
$1/d$ from lemma 2.2. Since ${\bf F}_{\mathfrak{p}}$ is a finite field of order $p^2$,  $|E_{\mathfrak{p}}({\bf
F}_{\mathfrak{p}})|= 1+p^2+a_p$ where $a_p=0, \pm p, \pm 2p$. The possibilities
$a_p = 0, \pm p$ arises only for $k={\bf Q}(i)$ and $k = {\bf Q}(\omega)$ 
respectively which we have omitted. Therefore we find that $|E_{\mathfrak{p}}(
{\bf F}_{\mathfrak{p}})| $ is either $ (p+1)^2$ or $(p-1)^2$. 
Hence  $p$ is congruent to either 1 or $-1$ modulo $M'$. The
density of such primes $p$ is $2/\varphi(M')$, and as before we get 
$$\frac{1}{d} \leq \frac{2}{\varphi(M')},$$
or $\varphi(M') \leq 2d$  completing the proof of the theorem.

\bigskip

\noindent {\bf Remark 2:} The technique of using supersingular primes can be
used in some situations to get better bounds for the order of the torsion
subgroup when $K$, the field of definition of $E$, is a non-normal
extension. For instance, if $L$ is a Galois extension of {\bf Q } with
Galois group $GL_2({\bf F}_q)$ which is disjoint from the field of CM,
$k$, and $K$ is the fixed field of the diagonal torus, then $[K:{\bf
Q}]=q(q+1)$. The set of primes $p$ in {\bf Q} with the property that there
is a prime ${\mathfrak{p}}$ in $K $ of degree 1 above $p$ corresponds to
those Frobenius substitutions in $GL_2({\bf F}_q)$ which have a conjugate
which is diagonalisable over ${\bf F}_q$. The set of elements in
$GL_2({\bf F}_q)$ which are diagonalisable over ${\bf F}_q$ can be easily
seen to be of cardinality $[(q-2)(q-1)q(q+1)]/2+(q-1)$. Therefore the set
of primes $p$ which are inert in $k$ and have a prime in $K$ of degree 1
above $p$ is roughly of density $1/4$. By the arguments in theorem 1.1 this
implies that the torsion in $E(K)$ is bounded by $M$ with $\varphi
(M)\leq 4$, implying that the set of possible values of $M$ is
$\{1,2,3,4,5,6,8\}$, instead of the much larger bound depending on $q$
coming from theorem 1.1. 

\end{section}

\begin{section}{A conjecture about torsion}

In this section we make a few general remarks and state a conjecture on
the bound for the order of a torsion point of an elliptic curve defined
over a number field. 

Let $X$ be a curve over {\bf Q} of genus $\geq 2$. Let $d$ be the gonality
of $X$, i.e., $d$ is the minimal integer among the degrees of maps from
$X$ to ${\bf P}^1$ and assume that $d$ is realised for a map $\pi$ defined
over {\bf Q}. It is clear that any element $x$ of $X(\overline{{\bf Q}})$
with $\pi (x)\in {\bf P}^1({\bf Q})$ is defined over a number field of
degree $\leq d$ and therefore there are infinitely many points in
$X(\overline{{\bf Q}})$ defined over number fields of degree $\leq d$.
Conversely, it has been proved by Debarre and Klassen in \cite{de-kl}
that if $X$ is a smooth {\it plane} curve then $d$ is the maximal integer
with the property that the number of points in $X(\overline{{\bf Q}})$
defined over a number field of degree $<d$ is finite. We would like to
believe that a suitably modified version of 
their theorem is valid also for modular curves $X_0(N)$ and
$X_1(N)$. More precisely, we believe that these modular curves have {\it no}
points, except for cusps,
 defined over an extension of {\bf Q} of degree $\leq Bd$, $B$ a
constant independent of $N$. 

It has been proved in \cite{Ab} that the gonality of $X_0(N)$ (resp.
$X_1(N)$) is at least a constant times the degree of the standard map of
$X_0(N)$ to ${\bf P}^1$; in fact, $d_0\geq (7/800)N$ (a similar bound but
quadratic in $N$ for $X_1(N)$). 
Earlier heuristics therefore lead us to the following.

\bigskip \noindent {\bf Conjecture:} {\em Let $E$ be an elliptic curve
defined over a number field $K$ with a torsion point of order $N$ in
$E(K)$. Then there is a constant $C$ independent of $E$ and $K$ such that
$\varphi (N)\leq C [K:{\bf Q}]$.  } \end{section}

\noindent {Mehta Research Institute, Chhatnag Road, Jhusi, Allahabad-211019, India
\\E-Mail: dprasad@mri.ernet.in}

\noindent{Olympiad Cell (NBHM), Dept. of Mathematics, Indian Institute of Science, 
Bangalore-560012, India\\ E-Mail:  yoga@math.iisc.ernet.in}

}
\end{document}